\documentstyle[12pt]{amsart}

\newtheorem{theorem}{Theorem}
\newtheorem{lemma}[theorem]{Lemma}

\newtheorem{algorithm}[theorem]{Algorithm}

\newcommand{\cc}{{\mathbb C}}

\newcommand\sA{{\cal A}}
\newcommand\sB{{\cal B}}

\newcommand\sE{{\cal E}}

\newcommand\sJ{{\cal J}}

\newcommand\sL{{\cal L}}

\newcommand\sO{{\cal O}}

\newcommand\sU{{\cal U}}

\newcommand\sX{{\cal X}}

\newcommand\sZ{{\cal Z}}

\newenvironment{proof}{\removelastskip\vskip12pt plus 1pt\noindent\em Proof.
\rm}{\hspace*{\fill}$\Box$\vskip12pt plus 1pt}

\begin{document}

\title[Numerical Homotopies to compute generic Points]
      {Numerical Homotopies to compute generic Points on positive
       dimensional Algebraic Sets}        

\author{Andrew J. Sommese}
\address{\hskip-\parindent
         Department of Mathematics \\
         University of Notre Dame \\
         Notre Dame, IN 46556, U.S.A.}

\email{sommese@@nd.edu and url: http://www.nd.edu/$\sim$sommese}

\author{Jan Verschelde}
\address{\hskip-\parindent
         Mathematical Sciences Research Institute \\
         1000 Centennial Drive \\
         Berkeley, CA 94720-5070, U.S.A.}
\address{\hskip-\parindent
         Department of Mathematics \\
         Michigan State University \\
         East Lansing, MI 48824-1027, U.S.A.}

\email{jan@@msri.org, jan@@math.msu.edu or jan.verschelde@@na-net.ornl.gov,
and urls: http://www.msri.org/people/members/jan
http://www.mth.msu.edu/$\sim$jan}

\thanks{Research at MSRI is supported in part by NSF grant DMS-9701755}

\date{June 28, 1999}

\begin{abstract}
\noindent Many applications modeled by polynomial systems have
positive dimensional solution components (e.g., the path synthesis
problems for four-bar mechanisms) that are challenging to compute
numerically by homotopy continuation methods. A procedure of A.
Sommese and C. Wampler
consists in slicing the components with linear subspaces in
general position to 
obtain generic points of the components as the isolated solutions of 
an auxiliary system.
Since this requires the solution of a number of larger
overdetermined systems, 
the procedure is computationally
expensive and also wasteful because many solution paths diverge.
In this article an embedding of the original polynomial system is
presented, which leads to a sequence of homotopies, with solution
paths leading to generic points of all components
as the isolated solutions of an auxiliary system. The new
procedure significantly reduces the number of paths to solutions
that need to be followed. This approach has been implemented and
applied to various polynomial systems, such as the cyclic n-roots
problem.

\bigskip

\noindent {\bf AMS Subject Classification.} 65H10, 68Q40.

\noindent {\bf Keywords}.  polynomial system, numerical homotopy
continuation, components of solutions, numerical algebraic
geometry, generic points, embedding.
\end{abstract}

\maketitle

\tableofcontents

\section{Introduction}

Let
\begin{equation}
 f(x) := \left[
  \begin{array}{c}
   f_1(x_1,\ldots,x_n) \\
   \vdots \\
   f_N(x_1,\ldots,x_n) \
 \end{array}\right]
\end{equation}
denote a system of $N$ polynomials on $\cc^n$.  Positive
dimensional components of the solution set of $f(x)=0$ are a
common occurrence, even when $N=n$. Sometimes they are an
unpleasant side show that happens with a system generated using a
model, for which only the isolated nonsingular solutions are of
interest; and sometimes,  the positive dimensional solution
components are of primary interest. In either case, dealing with
positive dimensional components, is usually computationally
difficult.  In \cite{SW}, Sommese and Wampler presented a
numerical  algorithm, which uses auxiliary systems to numerically
find  sets of solutions of the original system. These sets of
solutions, which let one numerically decide the dimension of the
zero set of the original system,  include at least the isolated
solutions of the original system, plus ``generic points'' of each
positive dimensional irreducible component of the solutions of the
original solutions. Generic points are the basic numerical data
which we  are using to investigate the positive dimensional
solution components.

\smallskip

The algorithm from \cite{SW}, which is based on
slicing with general linear spaces of different dimensions, leads
to $n$ auxiliary systems, which must be dealt with.  In this paper
we present an embedding of the system $f(x)=0$ into a family of
system of polynomials depending on $2n$ variables $(x,z)\in
\cc^{2n}$, and a large space of parameters.  We then single out
$n+1$ of the systems, $\sE_i(x,z)$ for $i$ from $n$ to $0$ (here
$\sE_0(x,z)=0$ is equivalent to the system $f(x)=0$) obtained by
choosing particular values of the parameters, plus a homotopy
$H_i$ going from $\sE_i$ to $\sE_{i-1}$. For simplicity we assume
that $f$ is not identically zero.  The system $\sE_n=0$ has
isolated nonsingular solutions. We use polynomial continuation
\cite{ag90,ag97}, \cite{li97}, \cite{mor87}
to implement the following algorithm

\begin{algorithm} {\rm Cascade of homotopies between embedded systems.

\smallskip

\begin{tabular}{lcr}
Input: $f$, $n$.
    & & {\em system with  solutions in $\cc^n$} \\
Output: $(\sE_i,\sX_i,\sZ_i)_{i=0}^n$.
    & & {\em embeddings with solutions} \\
\\
$\sE_0 := f$;
    & & {\em initialize embedding sequence} \\
for $i$ from 1 up to $n$ do
    & & {\em slice and embed} \\
\ \ \ $\sE_i$ := Embed($\sE_{i-1},z_{i}$);
    & & {\em $z_i$ = new added variable} \\
end for; & & {\em homotopy sequence starts} \\
$\sZ_n$ := Solve($\sE_n$);
    & & {\em all roots are isolated, nonsingular, with $z_n \not= 0$}\\
for $i$ from $n-1$ down to 0 do
    & & {\em countdown of dimensions} \\
\ \ \ $H_{i+1}$ := $
  t \sE_{i+1}
  +
  (1-t)\left(
     \begin{array}{c}
        \sE_i \\
        z_{i+1}
     \end{array}
  \right) $;
& &  {\em \begin{tabular}{r}
              homotopy continuation \\
       $t: 1 \rightarrow 0$ to remove $z_{i+1}$ \\
         \end{tabular} } \hspace{-5mm} \\
 \ \ \ $\sX_i$ := limits of solutions of $H_{i+1}$ \\
 \ \ \ \ \ \ as $t\to 0$ with $z_i=0$;
    & &  {\em on component} \\
\ \ \ $\sZ_i$ := $H_{i+1}(x,z_{i} \not= 0,t=0)$;
    & &  {\em not on component: these solutions} \\
    & &  {\em  are  isolated and nonsingular} \\
end for. & & \\
\end{tabular}
}
\end{algorithm}
The routine `Embed' will be defined in the next section.
In section three we present a worked out example of the algorithm.
Section four contains the mathematical background needed to prove
our main results:
\begin{enumerate}
\item
 if $i$ is the largest integer with $\sX_i$ nonempty, then
   the dimension of $f^{-1}(0)$ is $i$; and
\item  given any irreducible
    component $W$ of $f^{-1}(0)$ of dimension $i$, then, counting
    multiplicities,
   $\sX_i$ contains $\deg(W)$ generic points of $W$.
\end{enumerate}

\noindent The applications described in section five illustrate the
performance of the new procedure.
We end this paper with some directions for future research.

\bigskip

\noindent {\bf Acknowledgements.} Both authors would like to thank
MSRI for providing the stimulating environment where our
collaboration began.  We would also like to thank the Center for
Applied Mathematics, and the Duncan Chair of the University of
Notre Dame for their support of our work.  
The authors are grateful to G\"oran Bj\"orck for pointing at
the results of Uffe Haagerup.
The second author is supported by a post-doctoral fellowship at MSRI,
and supported in part by the NSF under Grant DMS - 9804846 at the
Department of Mathematics, Michigan State University.

\section{An Embedding of a Polynomial System} 
Throughout this paper we work with
 algebraic functions. This allows the possibility of using rational
  functions, and not
 just polynomials. This extra flexibility will be needed
 in a sequel, where, even though we start with a system of polynomials on
 $\cc^n$, it becomes necessary to work with rational functions on a
 Zariski open set of an associated Euclidean space.  We assume in
 what follows that locally we have the same number of equations as
 unknowns.  A procedure for reducing to this case is presented in
 \cite{SW}.

 We refer to \cite{SW} for a discussion of generic points.

 Given a system of algebraic functions
\begin{equation}
   f(x):=
   \left[
   \begin{array}{c}
   f_1(x)\\
   \vdots\\
   f_n(x)
   \end{array}
   \right]
\end{equation}
  on a connected algebraic manifold $X$ of dimension $n$
  embedded into a Euclidean space $\cc^A$
   we have the following basic embedding into a
  family of systems.  First we restrict
   $n$ linear functions on $\cc^A$ to $X$.
  Thus for $j$ from $1$ to $n$ we have
\begin{equation}
  L_j(x):=a_j+a_{j,1}x_1+\cdots+a_{j,A}x_A,
\end{equation}
  where $x_i$ is the restriction of the $i$-th coordinate function
  of $\cc^A$ to $X$.  By abuse of notation we let $L_j\in \cc^{A+1}$ denote
\begin{equation}
  \left(
  \begin{array}{cccc}
  a_j&a_{j,1}&\cdots&a_{j,A}\\
  \end{array}
  \right).
\end{equation}

  We fix linear coordinates $z_1,\ldots, z_n$ on a complex
  Euclidean space $\cc^n$. Then we have the system of equations
\begin{equation}
    \sE_i(f)(x,z,\Lambda_1,\ldots,\Lambda_i,L_1,\ldots,L_i):=
   \left[
    \begin{array}{c}
    f_1(x)+\sum_{j=1}^i\lambda_{1,j}z_j\\
    \vdots\\
    f_n(x)+\sum_{j=1}^i\lambda_{n,j}z_j\\
    a_1+a_{1,1}x_1+\cdots+a_{1,A}x_A+z_1\\
    \vdots\\
    a_i+a_{i,1}x_1+\cdots+a_{i,A}x_A+z_i
   \end{array}
   \right]
\end{equation}
where we let
\begin{equation}
  \Lambda_i :=
    \left[\begin{array}{c}
    \lambda_{1,i}\\
    \vdots\\
   \lambda_{n,i}
   \end{array}\right].
\end{equation}
We often refer to
$\sE(f)(x,z_1,\ldots,z_n,\Lambda_1,\ldots,\Lambda_n,L_1,\ldots,L_n)$
by $\sE_n$ or $\sE_n(f)$. Further we let
$\sE_i$ or
$\sE_i(f)$ denote
$\sE(f)(x,z_1,\ldots,z_i,\Lambda_1,\ldots,\Lambda_i,L_1,\ldots,L_i)$
on $\cc^{n+i}$.  Note that $\sE_0$ is just $f$; and that
  the solutions $(x,z_1,\ldots,z_i)\in X\times\cc^{i}$ of
\begin{equation}
  \sE_i(x,z_1,\ldots,z_i,\Lambda_1,\ldots,\Lambda_i,L_1,\ldots,L_i)=0
\end{equation}
 are naturally identified with the solutions
 $(x,z_1,\ldots,z_i,0\ldots,0)\in X\times\cc^{n}$ of the system
  \begin{equation}\label{identification}
 \left\{
   \begin{array}{ccc}
   \sE_i(x,z_1,\ldots,z_n,\Lambda_1,\ldots,\Lambda_i,0,\ldots,0,L_1,\ldots,L_i,0,\ldots,0)&=&0\\
   z_{i+1}&=&0\\
   \vdots&\vdots&\\
   z_n&=&0.
   \end{array}
 \right.
  \end{equation}

We let $Y$ denote the space $\cc^{n\times (A+1)}\times
\cc^{n\times n}$ of parameters
\begin{equation}
    \left[
   \begin{array}{ccccccc}
    a_1&a_{1,1}&\cdots&a_{1,A}&\lambda_{1,1}&\cdots&\lambda_{n,1}\\
     \vdots & \vdots & \ddots & \vdots & \vdots & \ddots & \vdots \\
    a_n&a_{n,1}&\cdots&a_{n,A}&\lambda_{1,n}&\cdots&\lambda_{n,n}\\
   \end{array}
   \right] \in \cc^{n\times (A+1)}\times\cc^{n\times n}
\end{equation}
for these systems. We have used the transpose of the
$\lambda_{i,j}$ for convenience in describing the stratification
$Y_0\subset Y_1\subset\cdots\subset Y_n$ of the space $Y$ given by
defining $Y_n:=Y$; and $Y_i$, for $i=0,\ldots,n-1$, as the subset
of $Y$ with the coordinates
\begin{equation}
    \left[
   \begin{array}{ccccccc}
    a_{i+1}&a_{i+1,1}&\cdots&a_{i+1,A}&\lambda_{1,i+1}&\cdots&\lambda_{n,i+1}\\
    \vdots & \vdots & \ddots & \vdots & \vdots & \ddots & \vdots \\
   a_n&a_{n,1}&\cdots&a_{n,A}&\lambda_{1,n}&\cdots&\lambda_{n,n}\\
   \end{array}
   \right] \in \cc^{(n-i)\times (A+1)}\times\cc^{(n-i)\times n}
\end{equation}
set equal to $0$.  Thus using the identification
(\ref{identification}) above, we can regard $Y_i$ as the parameter
space of the systems of equations
\begin{equation}
  \sE_i(x,z_1,\ldots,z_i,\Lambda_1,\ldots,\Lambda_i,L_1,\ldots,L_i)=0.
\end{equation}

 We consider homotopies $H_i$, for $i$ from $n$ to $1$ and $t$ from $1$ to
$0$,  defined by:
\begin{equation}
   H_i(x,z_1,\ldots,z_i,t) :=
    \left[
    \begin{array}{c}
    f_1(x)+\sum_{j=1}^{i-1}\lambda_{1,j}z_j+t\lambda_{1,i}z_i\\
    \vdots\\
    f_n(x)+\sum_{j=1}^{i-1}\lambda_{n,j}z_j+t\lambda_{n,i}z_i\\
    L_1(x)+z_1\\
    \vdots\\
    L_{i-1}(x)+z_{i-1}\\
    tL_i(x)+z_i
    \end{array}
   \right],
\end{equation}
with the convention that if $i=1$, we mean
\begin{equation}
  H_1(x,z_1,t) :=
  \left[
    \begin{array}{c}
    f_1(x)+t\lambda_{1,1}z_1\\
    \vdots\\
    f_n(x)+t\lambda_{n,1}z_1\\
    tL_1(x)+z_1\\
    \end{array}
   \right].
\end{equation}
Thus $H_i(x,z_1,\ldots,z_i,1)=\sE_i$ and
$H_i(x,z_1,\ldots,z_i,0)=0$ is
\begin{equation}
\left\{
  \begin{array}{ccc}
  \sE_{i-1}(x,z_1,\ldots,z_{i-1})&=&0\\
  z_i&=&0\\
  \end{array}
\right.
\end{equation}
Note that using this convention, $H_i$ can be rewritten as
$t\sE_i+(1-t)
 \left( \begin{array}{c} \sE_{i-1} \\ z_i \end{array} \right)$.

\begin{lemma}\label{keyLemma}
 There is a nonempty Zariski open set $U$ of points
\begin{equation}
   \left[
   \begin{array}{ccccccc}
    a_1&a_{1,1}&\cdots&a_{1,A}&\lambda_{1,1}&\cdots&\lambda_{n,1}\\
    \vdots & \vdots & \ddots & \vdots & \vdots & \ddots & \vdots \\
   a_n&a_{n,1}&\cdots&a_{n,A}&\lambda_{1,n}&\cdots&\lambda_{n,n}\\
   \end{array}
   \right] \in \cc^{n\times (A+1)}\times\cc^{n\times n}
\end{equation}
  such that for each $i=1,\ldots,n$,
 \begin{enumerate}
  \item the solutions of $\displaystyle \sE_i(x,z_1,\ldots,z_i)=0$ with
     $(z_1,\ldots,z_i)\not=0$ are isolated and nonsingular;
  \item given any irreducible
    component $W$ of $f^{-1}(0)$ of dimension $i$, and counting
    multiplicities, the isolated solutions of $\displaystyle
    \sE_i(x,z_1,\ldots,z_i)=0$ with $(z_1,\ldots,z_i)=0$,
    contain $\deg(W)$ generic points of $W$; and
  \item the solutions  of $\displaystyle \sE_i(x,z_1,\ldots,z_i)=0$ with
     $(z_1,\ldots,z_i)\not=0$ are the same as the solutions of
      $\displaystyle \sE_i(x,z_1,\ldots,z_i)=0$ with
     $z_i\not=0$.
 \end{enumerate}
\end{lemma}

\begin{proof}
 For a fixed $i\in\{0,\ldots,n\}$ consider the system of equations
\begin{equation}
\left\{
  \begin{array}{ccccc}
   \sum_{j=1}^n\alpha_{1,j}f_j & + & \sum_{j=1}^i\beta_{1,j}z_j &=&0\\
   &\vdots&\\
   \sum_{j=1}^n\alpha_{n,j}f_j & + & \sum_{j=1}^i\beta_{n,j}z_j &=&0 \\
   \gamma_1+\sum_{j=1}^A\gamma_{1,j}x_j & + & \sum_{j=1}^i\delta_{1,j}z_j  &=&0\\
   &\vdots&\\
   \gamma_i+\sum_{j=1}^A\gamma_{i,j}x_j & + & \sum_{j=1}^i\delta_{i,j}z_j &=&0. \
    \end{array}
\right.
\end{equation}
By Bertini's theorem (see section four for a convenient form of
this result), there is a nonempty Zariski open set $U$ of
\begin{equation}
  (\alpha,\beta,\gamma,\delta)\in \cc^{n\times n}\times \cc^{n\times
  i} \times \cc^{i\times (A+1)}\times \cc^{i\times i}
\end{equation}
  where
\begin{equation}
  \alpha:= \left[
  \begin{array}{ccc}
    \alpha_{1,1} & \cdots & \alpha_{1,n} \\
     \vdots & \ddots &  \vdots \\
    \alpha_{n,1} & \cdots & \alpha_{n,n}  \
  \end{array}
  \right]\ ;\ \
  \beta:= \left[
  \begin{array}{ccc}
    \beta_{1,1} & \cdots & \beta_{1,i} \\
     \vdots & \ddots & \vdots \\
    \beta_{n,1} & \cdots & \beta_{n,i}  \
  \end{array}
  \right]
\end{equation}
and;
\begin{equation}
  \gamma:= \left[
  \begin{array}{cccc}
   \gamma_1& \gamma_{1,1} & \cdots & \gamma_{1,A} \\
     \vdots & \vdots & \ddots & \vdots \\
   \gamma_i& \gamma_{1,1} & \cdots & \gamma_{i,A} \
  \end{array}
  \right]\ ;\ \
  \delta:= \left[
  \begin{array}{ccc}
    \delta_{1,1} & \cdots & \delta_{1,i} \\
     \vdots & \ddots & \vdots  \\
    \delta_{i,1} & \cdots & \delta_{i,i}  \
  \end{array}
  \right].
\end{equation}
  such that, if not empty, the zero set $Z$ of
the above system of equations on $X\times\cc^i$ minus the set of
common zeroes of $\displaystyle f_1,\ldots,f_n,z_1,\ldots,z_i$ is
smooth and of dimension $0$. Left multiplying this $n+i$ vector of
equations with an invertible $(n+i)\times (n+i)$ matrix $G$ of the
form
\begin{equation}
 \left[\begin{array}{cc}
   \sA & 0 \\
   0 & \sB \
 \end{array}
 \right]
\end{equation}
where $\sA$ is an $n\times n$ matrix and $\sB$ is a $i\times i$
matrix, results in the equivalent system of equations
\begin{equation}
  \begin{array}{ccccc}
   \sA\cdot\alpha \cdot f & + & \sA\cdot\beta\cdot\left[
                                             \begin{array}{c}
                                              z_1\\
                                              \vdots\\
                                              z_i
                                             \end{array}
                                             \right] &=&0\\
   \sB\cdot \gamma \cdot\left[
                         \begin{array}{c}
                          1\\
                          x_1\\
                         \vdots\\
                          x_A
                          \end{array}
                         \right] & + & \sB\cdot\delta\cdot\left[
                                             \begin{array}{c}
                                              z_1\\
                                              \vdots\\
                                              z_i
                                             \end{array}\right] &=&0.
     \end{array}
\end{equation}

Thus we can assume that the set $V$ is invariant under this action
by the matrices $G$.  Thus we have a nonempty Zariski open set $U$
of
\begin{equation}
   \left[
   \begin{array}{ccccccc}
    a_1&a_{1,1}&\cdots&a_{1,A}&\lambda_{1,1}&\cdots&\lambda_{n,1}\\
    \vdots & \vdots & \ddots & \vdots & \vdots & \ddots & \vdots \\
   a_n&a_{n,1}&\cdots&a_{n,A}&\lambda_{1,n}&\cdots&\lambda_{n,n}\\
   \end{array}
   \right] \in \cc^{n\times (A+1)}\times\cc^{n\times n}
\end{equation}
such that for each $i$ the equivalent system is of the form
\begin{equation}
 \left\{
    \begin{array}{ccc}
    f_1(x)+\sum_{i=1}^i\lambda_{1,i}z_i &=&0\\
    \vdots&\vdots\\
    f_n(x)+\sum_{i=1}^i\lambda_{n,i}z_i &=&0\\
    a_1+\sum_{j=1}^Aa_{1,j}x_j+z_1 &=&0\\
    \vdots&\vdots\\
    a_i+\sum_{j=1}^Aa_{i,j}x_j+z_i &=&0
   \end{array}
  \right.,
\end{equation}
and has smooth nonsingular zeroes when $(z_1,\ldots,z_i)\not=0$.
Thus we have the first assertion of the Lemma.

To see the second assertion of the Lemma, note that the solutions
of $\sE_i=0$ with $(z_1,\ldots,z_i)=0$ are naturally identified
with the solutions of the system
\begin{equation}
  \left\{
    \begin{array}{rcc}
    f_1(x)&=&0\\
    \vdots&&\vdots\\
    f_n(x) &=&0\\
    a_1+\sum_{j=1}^Aa_{1,j}x_j &=&0\\
    \vdots&&\vdots\\
    a_i+\sum_{j=1}^Aa_{i,j}x_j &=&0
   \end{array}
  \right.
\end{equation}
The second assertion follows now from the Algorithm in \cite[\S
3.1]{SW}.

We prove the third assertion of the lemma by induction on $i$. If
$i=1$, then  it is a tautology that the solutions  of
$\displaystyle \sE_i(x,z_1,\ldots,z_i)=0$ with
     $(z_1,\ldots,z_i)\not=0$ are the same as the solutions of
      $\displaystyle \sE_i(x,z_1,\ldots,z_i)=0$ with
     $z_i\not=0$.

So we can assume that the result is true for $k<i$ where $i>1$.
Note that a solution of $\sE_i=0$ with $z_i=0$ but
$(z_1,\ldots,z_i)\not = 0$ is a solution of the system
\begin{equation}
 \left\{
 \begin{array}{rcl}
 \sE_{i-1}(x,z_1,\ldots,z_{i-1})&=&0\\
 L_i(x)&=&0\\
 z_i&=&0\\
 &\vdots\\
 z_n&=&0.
 \end{array}
 \right.
\end{equation}
  Since the solutions of $\sE_{i-1}=0$ with $z_{i-1}\not=0$ are
isolated and nonsingular, a generic choice of $L_i(x)$ will not be
zero on any of the solutions. But this means that the solutions of
$\sE_i=0$  for generic $L_i(x)$ will have no solutions with
$z_i=0$.
\end{proof}

Note that if we choose $y$ generically  in $Y$, then we have
chosen the associated $y_i$ generically in $Y_i$ for each $i$ from
$n$ to $1$.  Thus we can assume that with an initial generic
choice of parameters $y$,  the behavior for each of the systems
$\sE_i=0$ is the behavior we expect with a generic choice of
parameters on $Y_i$.  One minor point remains.  Given a generic
choice of parameters $y\in Y$, generic behavior might not occur
for the homotopy $H_i$ with $t\in (0,1]$.  This is easily dealt
with by a trick of Morgan and Sommese \cite{MS1,MS2,MS3}.

Assume that we have chosen the parameters $y$ for the
systems in the nonempty Zariski open set $U$ of
\begin{equation}
   \left[
   \begin{array}{ccccccc}
    a_1&a_{1,1}&\cdots&a_{1,A}&\lambda_{1,1}&\cdots&\lambda_{n,1}\\
    \vdots & \vdots & \ddots & \vdots & \vdots & \ddots & \vdots \\
   a_n&a_{n,1}&\cdots&a_{n,A}&\lambda_{1,n}&\cdots&\lambda_{n,n}\\
   \end{array}
   \right] \in \cc^{n\times (A+1)}\times\cc^{n\times n}.
\end{equation}
We want our homotopies $H_i(x,z_1,\ldots,z_i,t)=0$ to define
algebraic sets that are flat over a Zariski open set containing
$(0,1]$.  Numerically this means that we want to  have
generic behavior for each $t\in (0,1]$, i.e., the same number of
 isolated points, no components that do not correspond to
 components in any other fiber. Exploiting generic flatness in this way
 is the underlying
 approach of the work of Morgan and Sommese \cite{MS1,MS2,MS3},
 to which we refer for more details.
 Since for all but
 a finite number of values of $t\in \cc$ generic behavior occurs, we can
 conclude that given a set of parameters $y\in U$ as above, the
 homotopy
\begin{equation}
  H_{\eta,i}(x,z_1,\ldots,z_i,t):= (x,z_1,\ldots,z_i,\eta t)=
    \left[
    \begin{array}{c}
    f_1(x)+\sum_{j=1}^{i-1}\lambda_{1,j}z_j+t\eta\lambda_{1,i}z_i\\
    \vdots\\
    f_n(x)+\sum_{j=1}^{i-1}\lambda_{n,j}z_j+t\eta\lambda_{n,i}z_i\\
    L_1(x)+z_1\\
    \vdots\\
    L_{i-1}(x)+z_{i-1}\\
    t\eta L_i(x)+z_i
    \end{array}
   \right],
\end{equation}
defines an algebraic set flat over a Zariski open set of $\cc$
containing $(0,1]$, for all but a finite set of $\eta\in \cc$ with
$|\eta|=1$.  Using openness of flatness, we can  absorb the $\eta$
into the parameters we use.  We still have a dense open set of
general parameters, but the set is only Zariski open in the
underlying real algebraic structure.

By Lemma~\ref{keyLemma} if the $L_i, \Lambda_i$ are chosen
randomly, then with probability one, $\displaystyle \sE_n(x,z)=0$
has a solution with $z=0$ only if $f(x)$ is identically zero on
$\cc^n$.

Recall from the introduction 
that for $i$ from $n$ to $1$,
\begin{enumerate}
\item $\sZ_i$ denotes the solutions to $\sE_i=0$ with $z_i\not=0$; and
\item  $\sX_{i-1}$ denote the limits with $z_{i-1}=0$
 of the paths of the
homotopy $H_i(t)$, from $t=1$ to $t=0$, starting at points of $\sZ_i$.
By convention, the condition $z_{i-1}=0$ is empty when $i=1$.
\end{enumerate} 

 \begin{theorem}\label{mainTheorem}
  Let $f$ be as above.  Assume that $f$ is not identically zero
  and that the $ \Lambda_i$ are chosen generically.
  If $i$ is the largest integer with $\sX_i$ nonempty, then
   the dimension of $f^{-1}(0)$ is $i$.  Moreover given any irreducible
    component $W$ of $f^{-1}(0)$ of dimension $i$, then, counting
    multiplicities,
   $\sX_i$ contains $\deg(W)$ generic points of $W$.
 \end{theorem}

Thus our algorithm achieves the same numerical goal of the
algorithm of \cite{SW}, but much more efficiently.

\noindent{\em Proof of Theorem \ref{mainTheorem}.}  We use the
notation of the proof of Lemma~\ref{keyLemma}.  By Lemma
\ref{keyLemma}, it follows that there is  a nonempty Zariski open
set of points $y\in Y$, such that for  each $i$ from $n$ to $1$,
all elements of the
 set $\sZ_i$ of solutions of $\sE_i(x,z_1,\ldots,z_i)=0$
with $z_i\not=0$ are nonsingular and isolated. Thus for a dense
set (Zariski open in the underlying real algebraic structure of
$y\in Y$)  the paths over $t\in (0,1]$ are smooth if they start at
nonsingular isolated solutions of $\sE_i=0$. By Lemma
\ref{localExtensionLemma}  and our random choice of $y\in Y$, each
nonsingular isolated solutions $(x^*,z_1^*,\ldots,z_{i-1}^*,0)$ of
$\sE_{i-1}=0$ must start from a nonsingular isolated solution
$(x',z_1',\ldots,z_i')$ of $\sE_i=0$.

 If $z'_i=0$ for a solution, then  $(z_1',\ldots,z_i')=0$ by
Lemma \ref{keyLemma}.
Thus $f(x')=0$ and hence $H_i(x',0,\ldots,0,t)=0$ identically in
$t$, and hence $(x^*,z_1^*,\ldots,z_{i-1}^*,0)=(x',0,\ldots,0)$.
Thus we know from Lemma \ref{keyLemma} that for a Zariski open set
of $Y_i$ the solutions of the system $\sE_i=0$ with $z_i=0$ (and
hence $(z_1,\ldots,z_i)=0$) are on irreducible components of the
algebraic set $f(x)=0$ of dimension $i$.  Thus
$(x^*,z_1^*,\ldots,z_{i-1}^*,0)=(x',0,\ldots,0)$ is on an
irreducible component of the algebraic set of points with
$\sE_{i-1}=0$, which has dimension $\ge 1$.  This contradicts
 $(x^*,z_1^*,\ldots,z_{i-1}^*)$ being an isolated nonsingular solution
 of $\sE_{i-1}=0$.  Thus the isolated nonsingular solutions of
 $\sE_{i-1}=0$ must be endpoints of paths of homotopy $H_i$
 starting at points of $\sZ_i$.

 Using Lemma \ref{localExtensionLemma} the above argument
 can be modified  to show that isolated, but possibly singular,
 solutions of  $(x^*,z_1^*,\ldots,z_{i-1}^*,0)$ of $\sE_{i-1}=0$ must
start from a nonsingular isolated solution $(x',z_1',\ldots,z_i')$
of $\sE_i=0$ with $(z_1',\ldots,z_i')\not=0 $ (and hence by Lemma
\ref{keyLemma} with $z_i'\not=0$). \qed

 The property of a generic system of the form $\sE_n(f)$ that
it has isolated nonsingular solutions is useful.  In this
direction, see \cite{lsy89}.

\section{A Worked Out Example}

Consider the polynomial system $f(x) = 0$ and start system $g(x) = 0$.
\begin{equation} \label{eqextarstar}
  f(x) =
  \left\{
     \begin{array}{r}
        x_1^2 x_2 = 0 \\
        x_1^2 ( x_2^2 + x_1 ) = 0
     \end{array}
  \right.
  \quad \quad
  g(x) =
  \left\{
     \begin{array}{r}
        x_1^3 - c_1 = 0 \\
        x_2^4 - c_2 = 0
     \end{array}
  \right.
\end{equation}
To solve $f(x) = 0$ we trace $D = 3 \times 4 = 12$
solution paths starting at the solutions of $g(x) = 0$.
One path diverges to infinity, three paths converge to
$(0,0)$, and the remaining eight paths end at the solution
component $x_1 = 0$, for some $x_2 \not= 0$.  The condition
numbers at the end of the paths do not allow us to decide which
solution is isolated.

\smallskip

To embed $f(x) = 0$ we take a random hyperplane
$L(x) = a_0 + a_1 x_1 + a_2 x_2 = 0$
and choose two random complex constants,
$\lambda_1$ and $\lambda_2$:
\begin{equation} \label{eqembedsys}
  \sE_1(x,z) =
  \left\{
     \begin{array}{r}
         x_1^2 x_2 + \lambda_1 z = 0 \\
         x_1^2 ( x_2^2 + x_1 ) + \lambda_2 z = 0 \\
         a_0 + a_1 x_1 + a_2 x_2 + z = 0
     \end{array}
  \right.
\end{equation}
We can solve this system by tracing $D = 12$ solution paths, using
a standard linear homotopy with the equations $g(x) = 0$, $z - 1 = 0$
as start system, with $g(x) = 0$ as in~(\ref{eqextarstar}).
Five paths diverge to infinity. Two
paths go to the same solution with $z = 0$, which reveals the
degree of the solution component~$x_1^2 = 0$.  Note that
geometrically this component corresponds to a pair of lines. The
five remaining paths go to regular solutions with $z \not= 0$.

\smallskip

To compute the possible remaining isolated solutions, we trace
five solution paths starting at the five regular solutions of the
system~(\ref{eqembedsys}).  In going with $t$ from 1 to 0,
we use the homotopy

\begin{equation} \label{eqembedhom}
  H_1(x,z,t) = t \left(
  \left\{
    \begin{array}{r}
       x_1^2 x_2 + \lambda_1 z = 0 \\
       x_1^2 ( x_2^2 + x_1 ) + \lambda_2 z = 0 \\
       a_0 + a_1 x_1 + a_2 x_2 + z = 0
    \end{array}
  \right.
  \right)
  + (1-t)
  \left(
  \left\{
     \begin{array}{r}
        x_1^2 x_2 = 0 \\
        x_1^2 ( x_2^2 + x_1 ) = 0 \\
        z = 0
     \end{array}
  \right.
  \right)
\end{equation}
Three paths converge to $(0,0)$, two paths go to solutions on the
component with $x_1 = 0$ and $x_2 \not= 0$. End games are still
needed to decide whether the solutions are isolated.

\smallskip

With our new method, 17 solution paths instead of 24 were
traced, as 24 would have been the number of paths with an
iteration of the procedure in~\cite{SW}.

\smallskip

Polyhedral root counting methods provide a generically sharp root
count for polynomial systems.  In particular, the mixed volume of
the Newton polytopes of the system equals the number of roots in
$(\cc^*)^n$, $\cc^* := \cc \setminus \{ 0 \}$, for a system with
generic coefficients.  When the system has only few monomials with
nonzero coefficients, then the mixed volume provides a much lower
root count than the B\'ezout bounds based on the degrees of the
polynomials.

\smallskip

For our type of applications, the distinction between solutions with
$z = 0$ and $z \not= 0$ is instrumental in identifying components of
solutions.
This difference does depend on the values of the
coefficients of the original system and is not neglected by the
ordinary mixed volume.
So we will not miss any solutions with $z = 0$, but we may miss
solution components for which some $x_i = 0$.
Fortunately, extensions of the polyhedral methods that
allow to count and compute all affine roots (that is in $\cc^n$
instead of $(\cc^*)^n$) are covered amply in the literature
(see \cite{ev97}, \cite{glw98}, \cite{hs97}, \cite{lw96},
\cite{roj94,roj99} and \cite{rw96}).
The key idea~\cite{lw96} is to add a random constant to every equation
to shift the roots with zero components away from the coordinate axes.
In removing these constants by continuation, all affine roots lie at
the end of some path that starts at a root in $(\cc^*)^n$.

\smallskip

Consider again the polynomial system $f(x) = 0$ in~(\ref{eqextarstar}).
Because of the first equation we immediately see that there cannot be
any solution with all components different from zero.
The direct application of polyhedral methods does not yield anything,
since the mixed volume for $f(x) = 0$ equals zero.
With affine polyhedral methods, we consider the system
\begin{equation}
   f^{(0)}(x) =
   \left\{
      \begin{array}{r}
         x_1^2 x_2 + \gamma_{1} = 0 \\
         x_1^2 ( x_2^2 + x_1 ) + \gamma_{2} = 0
      \end{array}
   \right.
\end{equation}
where $\gamma_{1}$ and $\gamma_{2}$ are randomly chosen complex
numbers.  Here the mixed volume equals five.  Note the difference
with the total degree $D = 12$.
Letting the $\gamma$'s go to zero, three of the five paths converge
to the origin, and the other two paths go to other solutions on the
component.

\smallskip

To obtain information about the components, with polyhedral methods
we consider the embedding
\begin{equation}
   \sE_1(x,z) =
   \left\{
      \begin{array}{r}
         x_1^2 x_2 + \gamma_{1} + \lambda_1 z = 0 \\
         x_1^2 ( x_2^2 + x_1 ) + \gamma_{2} + \lambda_2 z = 0 \\
         a_0 + a_1 x_1 + a_2 x_2 + z = 0
      \end{array}
   \right.
\end{equation}
This system has mixed volume equal to seven.  We use this as start
system to solve it with $\gamma_{1} = 0$ and $\gamma_{2} = 0$,
following the paths that start at the seven solution paths of
$\sE_1(x,z) = 0$.  From the seven paths, five paths
go to solutions with $z \not= 0$, and the other two paths go to
the component ending with $z = 0$.

\section{Bertini's Theorem and a Local Extension Theorem}

Here is a weak, but convenient form of Bertini's theorem,
e.g., Fulton \cite[Example 12.1.11]{Fulton}.  For a further
discussion of Bertini theorems, see also \cite[\S 1.7]{BeSo}.

\begin{theorem}[Bertini]\label{bertini}
Let $X$ be an algebraic manifold of dimension $n$, e.g., complex
Euclidean space or complex projective space. Let $\sL_1,\ldots,
\sL_n$ be line bundles on $X$.  For $i$ from $1$ to $n$, let
$\displaystyle\{s_{i,j}| j=1,\ldots,r_i\}$ be a set of sections of
$\sL_i$.  Let $B_i$ denote the set of common zeroes of  the
sections $\displaystyle\{s_{i,j} | 1\le j\le r_i\}$; and let
$B:=\cup_i B_i$. Then given general real or complex numbers
$\displaystyle\{\lambda_{i,j}| j=1,\ldots,r_i; i=1,\ldots,n\}$,
all solutions on $X-B$, of the system of equations
\begin{equation}
\left\{
  \begin{array}{ccc}
     {\displaystyle \sum_{j=1}^{r_1}\lambda_{1,j}s_{1,j}} & = & 0 \\
                                       & \vdots & \\
     {\displaystyle \sum_{j=1}^{r_n}\lambda_{n,j}s_{n,j}} & = & 0
 \end{array}
\right.
\end{equation}
are isolated and nonsingular.
\end{theorem}

In the proof of Theorem \ref{mainTheorem}, we need an extension
theorem, giving conditions ensuring that if we have a system of
equations $f(x,y)=0$ depending on parameters $y$ with a
multiplicity $\mu$ isolated solution $x^*$ of the system for some
point $y^*$ in the parameter space, then there is a neighborhood
$U$ of $x^*$ such that for  points $y$ near $y^*$ in the parameter
space, there are $\mu$ solutions counting multiplicities of the
system $f(x,y)=0$. Special cases are well known in the context of
all polynomial systems, but we do not know a general reference in
the numerical analysis literature.  This sort of result is
standard for algebraic geometers in the algebraic context, or
within the German school of several complex variables in the
complex analytic context.  So we are simply explaining why this
sort of result follows immediately from standard results in these
fields. We work locally.  All open sets are in the usual Euclidean
topology, i.e., not in the Zariski topology.

A possibly nonreduced complex analytic space, $\sX$, is said to be
a local complete intersection if given any point $x\in \sX$, there
is a set $\sU\subset \sX$, open in the usual complex topology,
that contains $x$, and such that
 \begin{enumerate}
 \item there is an embedding $\phi : \sU\to B$ of $\sU$ into an open ball
 in $B\subset
 \cc^N$ for some $N:=n+m>0$;
 \item the ideal of the complex space $\phi(\sU)$ is defined by
 holomorphic functions $g_1,\ldots,g_n$; and
 \item the dimension of the maximal dimensional irreducible
 component of $\sX$ through $x$ is $m$.
 \end{enumerate}
Local complete intersections are very well behaved.  One
elementary, but important fact about them is that all the
irreducible
 component of $\sX$ through a given $x\in \sX$ are equal dimensional.
 A less elementary, but important result  is that
 for surjective  holomorphic maps (respectively algebraic morphisms)
 from possibly
 nonreduced complex analytic spaces (respectively possibly nonreduced algebraic
 varieties) to complex manifolds (respectively algebraic
 manifolds), flatness of the morphism is equivalent to the
 openness of the morphism (which is equivalent to all the irreducible
 components of all
 the fibers of the
 morphism restricted to a given connected component of $\sX$ being equal
 dimensional).  This last result is a standard result for complex
 analytic spaces
  $\sX$, whose structure sheaf consists of local rings which are Cohen-Macaulay,
  e.g., Fischer \cite[page 158]{Fischer}.
 Local complete intersections are among the simplest examples after manifolds of
 complex analytic spaces with  Cohen-Macaulay
 local rings.

 Let $X$ be a connected $n$-dimensional complex manifold.  Let $Y$
be a connected $m$-dimensional complex manifold. Let
\begin{equation}
  f(x,y)=
  \left[
  \begin{array}{c}
   f_1(x,y)\\
    \vdots\\
    f_n(x,y)
  \end{array}
  \right]=0
\end{equation}
be a system of $n$ holomorphic functions.  Let $x^*$ be an
isolated solution of the system $f(x,y^*)=0$ for a fixed value
$y^*\in Y$, i.e., assume that there is an open set $\sO\subset X$
containing $x^*$ with $x^*$ the only solution of $f(x,y^*)=0$ on
$\sO$.  Assume that the multiplicity of $x^*$ is $\mu$.  This
number, which is $1$ exactly when the Jacobian of $f(x,y^*)$ is
invertible at $x^*$, is equal to
\begin{equation}
  \dim_\cc \sO_{X|x^*}/\sJ(f_1(x,y^*),\ldots,f_n(x,y^*))
\end{equation}
where $\sO_{X|x^*}$ is the local ring of convergent power series
on $X$ centered at the point $x^*$, and
$\sJ(f_1(x,y^*),\ldots,f_n(x,y^*))$ is the ideal in $\sO_{X|x^*}$
generated by the functions $f_1(x,y^*),\ldots,f_n(x,y^*)$.

\begin{lemma}[Local Extension Lemma]\label{localExtensionLemma}Let
$X$, $Y$, $f$, $x^*$, and  $y^*$ be as above. There are open
neighborhoods $U$ of $x^*\in X$ and $V$ of $y^*\in Y$ such that
for any $y\in V$ there exist $\mu$ isolated solutions (counting
multiplicities) of $f(x,y)=0$ on $U$.
\end{lemma}
 \begin{proof} Let $\dim_{\{(x^*,y^*)\}} Z$ denote the dimension of an
 analytic set
$Z\subset X\times Y$ at the point $(x^*,y^*)$. Let $X'$ denote the
zero set of $f_1,\ldots,f_n$ on $X\times Y$. Since there are
$n$-functions, the dimension of each irreducible component of $X'$
is at least $m$, and in particular
\begin{equation}
  \dim_{(x^*,y^*)} X'\ge m.
\end{equation}
   Since $x^*$ is isolated we
have that $\dim_{(x^*,y^*)} (X\times \{y^*\})\cap X'=0$. Since
$X\times Y$ is smooth, we have
 \begin{eqnarray*}
 \dim_{(x^*,,y^*)} (X\times \{y^*\})\cap X'&\ge& \dim_{(x^*,,y^*)}X\times
 \{y^*\})+\dim_{(x^*,,y^*)}X'-\dim_{(x^*,,y^*)}X\times Y\\
 &=&n+\dim_{(x^*,,y^*)}X'-n-m.
 \end{eqnarray*}
Thus we conclude that
 \begin{equation}\label{locCompleteInt}
 \dim_{(x^*,y^*)} X'= m.
 \end{equation}
Thus the union $X''$ of the components of $X'$ passing through
$(x^*,y^*)$ has pure dimension $m$. It follows, e.g., use Gunning
\cite[Theorem 16]{Gunning}, that there are neighborhoods $U$ of
$x^*\in X$ and $V$ of $y^*\in Y$ such that the projection $\pi$ of
$X''\cap (U\times V)$ to $V$ is proper and finite.  By shrinking
$U$ and $V$ we can assume that $X''\cap (U\times V)=X'\cap
(U\times V)$ and that  $X'\cap (U\times \{y^*\})=(x^*,y^*)$.

Now $X'\cap (U\times V)$ is a local complete intersection, and
thus since the map $\pi_{X'}: X'\to V$ is proper with finite
fibers, we conclude by the discussion before the discussion of
flatness, that this map is flat. Thus the direct image
$\sE:=\pi_*(\sO_{X'})$ of the structure sheaf of $X'$  is locally
free, e.g., Fischer \cite[Proposition 3.13]{Fischer}.  On $U\times
V$,  $X'$  is defined by the functions $f_i(x,y)=0$, for
$i=1,\ldots,n$.  By definition, at a point $(x',y')\in (U\times
V)\cap X'$ we have $\sO_{X'|\{(x',y')\}}$ is
\begin{equation}
  \sO_{(U\times V)|\{(x',y')\}}/\sJ(f_1(x,y),\ldots,f_n(x,y))
\end{equation}
 where $\sO_{(U\times V)|\{(x',y')\}}$ is the local ring of convergent power
series on $U\times V$ centered at $\{(x',y')\}$ and
$\sJ(f_1(x,y),\ldots,f_n(x,y))$ denotes the ideal in
$\sO_{(U\times V)|\{(x',y')\}}$ generated by $f_1,\ldots,f_n$.

The statement that $\sE$ is locally free, is equivalent to ranks
of $\sE$ at different points of $V$ being equal.   The rank of
$\sE$ at a point $y'\in V$ is by definition
\begin{equation}
  \dim_\cc\sE_{y'}:=
  \sE/\left({\mathfrak m}_{y'}\otimes_{\sO_{Y|y'}}\sE\right)
\end{equation}
  where ${\mathfrak m}_{y'}$ is the maximal ideal
generated $\sO_{Y|y'}$ consisting of convergent power series
vanishing at $y'$.  Thus comparing to the definition of the
multiplicity of the solution $x^*$ of $f(x,y^*)$ being $\mu$, we
see that the rank of $\sE$ is $\mu$.  Unwinding the definition of
$\pi_*$, $\sE_{y'}$ for any fixed $y'\in V$ is the direct sum of
the modules
\begin{equation}
  \sO_{X|x'}/\sJ(f_1(x,y'),\ldots,f_n(x,y')),
\end{equation}
with index set the set of distinct points $x'\in U$ with
$f(x',y')=0$.  This proves the assertion of the Lemma.
\end{proof}

\section{Applications and Computational Experiences}

We have conducted a systematic set of experiments with PHCpack~\cite{ver97}
on four case studies of familiar polynomial systems.
Little symbolic manipulation of polynomials is needed to set up homotopies
derived from the embedding presented in this paper.

\smallskip

To avoid wasting paper, we have omitted the algebraic formulations
of the systems, which can either be found electronically in the
database maintained on the web sites of the second author, or
can be consulted in the cited literature.
Although the systems are academic examples and only interesting
for benchmarking purposes, we try to indicate the relevance to
their application fields.

\smallskip

Unless stated otherwise, all reported timings concern a
450MHz Intel Pentium machine with 1Gb of main memory and 1Gb
swap space, running Debian GNU/Linux.
We observe that the initial stage is always the most expensive one.
Dealing with components of solutions is a much harder
problem than just approximating the isolated solutions.

\smallskip

With a faster computer one can attack larger problems, but also
the quality of the software matters.  Concerning this latter aspect
we want to point out that no special-purpose software has been
developed that exploits the structure of the embedded systems.
The continuation could for instance go faster if one eliminated
explicitly some variables using the linear hyperplanes that were
introduced in the slicing.

\subsection{A planar four-bar mechanism}

Four-bar mechanisms are ubiquitous in mechanical design.
The 4-variable polynomial system that was derived and solved in~\cite{mw90}
has total degree 256 and lowest multi-homogeneous B\'ezout bound 96.
The mixed volume equals 80, avoiding the calculation of zero component
solutions.  There is a solution component of dimension two, with
sum of degrees equal to two.
The lowest multi-homogeneous B\'ezout bound of $\sE_2$
equals 240 whereas the mixed volume equals~96.

\smallskip

We summarize the computational experiments in Table~\ref{tabfourbar},
whose format goes as follows.  The cascade of homotopies starts with
solving the start system $g$ by polyhedral methods since mixed volumes
are sharper than the bounds based on the degrees.
To solve the system $\sE_k$ we must trace \#paths solution curves.
This number is partitioned into ones on the component ($z=0$),
regular finite solutions ($z \not= 0$), and diverging ones
($\rightarrow \infty$).  Observe the take over of the $(z \not= 0)$-solutions
to the next row.
The number in the column with heading $z \not= 0$ on the $\sE_0$-row
is the number of isolated solutions of the original system.
In the last column we list cpu times.

\begin{table}[hbt]
\begin{center}
\begin{tabular}{|c||r|r|r|r||r|} \hline
system & \#paths & $z = 0$ & $z \not= 0$ & $\rightarrow \infty$
& cpu time~~ \\ \hline \hline
 $g$, start & ~96 & ~0 & ~96 & ~0  &   10s~880ms \\
 $\sE_2$    & ~96 & ~2 & ~68 & 26 & 4m~47s~830ms \\
 $\sE_1$    & ~68 & 12 & ~56 & ~0 &    18s~370ms \\
 $\sE_0$    & ~56 & 20 & ~36 & ~0 &    11s~~10ms \\ \hline \hline
  total     & 316 & 34 & 256 & 26 & 5m~28s~~90ms \\ \hline
\end{tabular}
\caption{The number of paths and timings for a planar four-bar mechanism.}
\label{tabfourbar}
\end{center}
\end{table}

The information in Table~\ref{tabfourbar} is useful to make some rough
comparisons with other homotopy methods.  For instance, if we are only
interested in the isolated roots, the cost would be of the same magnitude
as on the $\sE_0$-row.  To compare with the procedure of Sommese
and Wampler~\cite{SW}, we could add up the entries in the columns with
headings $z = 0$ and $\rightarrow \infty$ to the \#paths in the next
row, simulating the process of no solution recycling.

\smallskip

We must note that there exists an isotropic formulation of this
problem~\cite{wam96} for which the mixed volume is an exact root count
for the 36 isolated solutions.
The black-box solver of PHC needs only 11s~710ms to solve that system.

\subsection{Constructing Runge-Kutta formulas}

The following system is due to Butcher~\cite{but84} and was used
as test problem in~\cite{bgk86}.
The system has seven variables and a three-dimensional component.
Its total degree is 4608 and mixed volume equals 24.
After slicing and embedding thrice, the resulting 10-variable
system has mixed volume 247.
We summarize the running statistics in Table~\ref{tabbutcher}.

\begin{table}[hbt]
\begin{center}
\begin{tabular}{|c||r|r|r|r||r|} \hline
system & \#paths & $z = 0$ & $z \not= 0$ & $\rightarrow \infty$
& \multicolumn{1}{c|}{cpu time} \\ \hline \hline
 $g$, start & 247 & ~0 & 247 & ~~0 & 17m~13s~770ms\\
 $\sE_3$    & 247 & ~3 & 193 & ~51 & ~9m~57s~~60ms \\
 $\sE_2$    & 193 & 15 & 161 & ~17 & ~1m~59s~470ms \\
 $\sE_1$    & 161 & 11 & ~72 & ~78 & ~4m~57s~600ms \\
 $\sE_0$    & ~72 & ~0 & ~~4 & ~68 & ~2m~24s~650ms \\ \hline \hline
   total    & 920 & 29 & 677 & 214 & 34m~32s~550ms \\ \hline
\end{tabular}
\caption{The number of paths and timings for Butcher's system.}
\label{tabbutcher}
\end{center}
\end{table}

\subsection{Modular equations for special algebraic number fields}

The authors of~\cite{gt89} discuss the finding of the solution
components of the system that was derived in~\cite{cd88}.
That system has a component of dimension two, and in~\cite{cd88} the
authors reasoned that any three of the four polynomials generate the
same ideal.  The challenge posed in~\cite{cd88} to computer algebra
was to verify this claim.

\smallskip

Although we believe that nowadays, more than 10 years after~\cite{cd88},
some computer algebra packages might be strong enough to compute
a Gr\"obner basis for this problem,
we would like to indicate how one could answer the challenge numerically.
The approach we propose is to replace any of the four equations by
the same random hyperplane.  This operation cuts the dimension by one.
So we have to embed each of the four systems with the same random
constants.  We solve one of them, this takes about 2 minutes,
and check whether the generic points we find satisfy the other equations.

\subsection{On Fourier transforms: the cyclic n-roots problem}

In~\cite{bf91} an example was given which stems from the problem of
finding all ``bi-equimodular'' vectors ${\bf x} \in \cc^n$, i.e.:\
all ${\bf x}$ with coordinates of constant absolute value such that
the Fourier transform of $\bf x$ is a vector with coordinates
of constant absolute value, see also~\cite{bjo85} and \cite{bjo89}.
This system was popularized in~\cite{dav87}, and is by far the most
notorious benchmark problem in polynomial system solving.
Other references are~\cite{bf91a},~\cite{bf94} and~\cite{ec95}.

\smallskip

In~\cite{moe98}, the following conjectures of Ralf Fr\"oberg are mentioned.
If $n$ has a quadratic divisor, then there are infinitely many solutions.
If the number of solutions is finite, then this number equals
${\scriptsize \left( \begin{array}{c} 2n-2 \\ n-1 \end{array} \right)}$.
Uffe Haagerup~\cite{haa96} has proven that 
for $n$ prime, the number of roots is always finite, and equals
indeed $\frac{(2n-2)!}{(n-1)!^2}$.

\smallskip

In Table~\ref{tabcyclic8} we summarize the results of the embedding
for the cyclic 8-roots problem, which has a solution component of
dimension one.  The original system has mixed volume equal to 2,560.
The embedded system has mixed volume 4,176.
See Table~\ref{tabcyclic8} for the summary of the runs.
We remark that if one is only interested in the isolated roots,
one only has to trace 2,560 solution paths.

\begin{table}[hbt]
\begin{center}
\begin{tabular}{|c||r|r|r|r||r|} \hline
system & \#paths & $z = 0$ & $z \not= 0$ & $\rightarrow \infty$
& \multicolumn{1}{c|}{cpu time} \\ \hline \hline
 $g$, start  & ~4,176 & ~~0 & 4,176 & ~~~~0 & 1h~30m~26s~800ms \\
 $\sE_1$     & ~4,176 & 144 & 3,975 & ~~~57 & 1h~10m~~3s~930ms \\
 $\sE_0$     & ~3,975 & 495 & 1,152 & 2,328 & 7h~~7m~55s~580ms \\ \hline \hline
  total      & 12,327 & 639 & 9,303 & 2,385 & 9h~48m~26s~310ms \\ \hline
\end{tabular}
\caption{The number of paths and timings for cyclic 8-roots.}
\label{tabcyclic8}
\end{center}
\end{table}

Before developing the embedding method, we attacked the cyclic 8-roots
problem first with the slicing method of Sommese and Wampler.
The computation of generic points on the solution component required
the tracing of 10,940 solution paths (which takes now 4,176 paths).
For this task, the black-box solver of PHC needed 4 days and nights
15h~6m~49s~392ms cpu time on a SUN workstation.
Although the Linux machine we used is somewhat faster than that
SUN workstation, it certainly does not speed up things 50 times!

\smallskip

We end this section with another illustration on the importance of
having good polynomial equations for the same problem.
Thanks to a trick of John Canny, described in~\cite{emi94},
there exists a reduced version of the cyclic $n$-roots problem.
This version reduces the number of roots by eight and
yields a significant saving in the path tracking.
For instance, the mixed volume of the reduced cyclic 8-roots problem
equals~320, instead of~2,560 for the original problem.
Table~\ref{tabredcyc8} summarizes the computational experiments.

\begin{table}[hbt]
\begin{center}
\begin{tabular}{|c||r|r|r|r||r|} \hline
system & \#paths & $z = 0$ & $z \not= 0$ & $\rightarrow \infty$
& \multicolumn{1}{c|}{cpu time} \\ \hline \hline
 $g$, start  & ~~775 & ~~0 & ~~775 & ~~0 & 10m~~1s~860ms \\
 $\sE_1$     & ~~775 & ~18 & ~~744 & ~13 & 11m~23s~680ms \\
 $\sE_0$     & ~~744 & 445 & ~~144 & 155 & 14m~52s~340ms \\ \hline \hline
  total      & 2,294 & 463 & 1,663 & 168 & 36m~17s~860ms \\ \hline
\end{tabular}
\caption{The number of paths and timings for reduced cyclic 8-roots.}
\label{tabredcyc8}
\end{center}
\end{table} 

Finally, this reduction also allows us to treat the cyclic 9-roots
system that has a two-dimensional component of solutions.
The sum of the degrees of the reduced cyclic 9-roots problem equals two.
The mixed volume of embedded system $\sE_2$ equals 4,044, computed in
cpu time 2h~48m~4s~320ms.  The polyhedral continuation for the start
system requires 1h~4m~5s~130ms and the path following to solve $\sE_2$
takes 2h~15m~4s~540ms.  In total, the black-box solver of PHCpack
needs 6h~7m~45s~270ms to solve~$\sE_2$.

\section{Conclusions and Future Directions}

In the paper we have presented an embedding of polynomial systems
to compute generic points on components of solutions.
The homotopies recycle solutions while moving down to the isolated
solutions.  The major advantage compared to the method in~\cite{SW}
is that fewer solution paths diverge.  We provide practical evidence
for the efficiency of the embedding technique.

\smallskip

Besides the reported progress on homotopy methods, we want to point
out that our embedding technique might be of independent interest,
in the sense that it allows 
us to focus on one component of a particular
dimension.  The main novelty of the embedding consists in the fact
that the solutions of the embedded system that 
do not lie on the focused
component provide meaningful information 
about the components of lower
dimension.  We hope that besides homotopy methods, other solvers for
polynomial systems could benefit from our embedding technique.

\smallskip

It is sometimes said that every scientific problem solved
generates three other questions.  Here we briefly indicate three
different directions for future research.  First of all, the
embedding leads to higher-dimensional systems with special
structure.  It would be worthwhile to exploit this structure to
obtain a more efficient mixed-volume computation to alleviate the
complexity of the first stage. Secondly, the computed generic
points on the components are an interesting starting point for
future computational explorations of the solution components.
Note that continuation is a very natural tool to scan the
component while wiggling the added hyperplanes. Lastly, in the
final stage of the cascade of embedded homotopies, it remains
difficult to separate solution paths converging to components from
other possibly ill-conditioned isolated solutions. End games that
produce certificates in the form of a random hyperplane through a
nonisolated solution will have to be developed.


\begin{thebibliography}{10}

\bibitem{ag90}
E.L. Allgower and K.~Georg.
\newblock {\em Numerical {C}ontinuation {M}ethods, an {I}ntroduction},
  volume~13 of {\em Springer Ser. in Comput. Math.}
\newblock Springer--Verlag, Berlin Heidelberg New York, 1990.

\bibitem{ag97}
E.L. Allgower and K.~Georg.
\newblock Numerical {P}ath {F}ollowing.
\newblock In P.G. Ciarlet and J.L. Lions, editors, {\em Techniques of
  Scientific Computing (Part 2)}, volume~5 of {\em Handbook of Numerical
  Analysis}, pages 3--203. North-Holland, 1997.

\bibitem{bf91a}
J.~Backelin and R.~Fr\"{o}berg.
\newblock How we proved that there are exactly 924 cyclic 7-roots.
\newblock In {\em Proceedings of ISSAC-91}, pages 101--111. ACM, 1991.

\bibitem{BeSo}
M.~Beltrametti and A.J. Sommese.
\newblock {\em The adjunction theory of complex projective varieties},
  volume~16 of {\em Expositions in Mathematics}.
\newblock Walter De Gruyter, Berlin, 1995.

\bibitem{bjo85}
G.~Bj\"{o}rck.
\newblock Functions of modulus one on {$Z_p$} whose {F}ourier transforms have
  constant modulus.
\newblock In {\em Proceedings of the Alfred Haar Memorial Conference,
  Budapest}, volume~49 of {\em Colloquia Mathematica Societatis J\'anos
  Bolyai}, pages 193--197. 1985.

\bibitem{bjo89}
G.~Bj\"{o}rck.
\newblock Functions of modulus one on {$Z_n$} whose {F}ourier transforms have
  constant modulus, and ``cyclic n-roots''.
\newblock In J.S. Byrnes and J.F. Byrnes, editors, {\em Recent Advances in
  Fourier Analysis and its Applications}, volume 315 of {\em NATO Adv. Sci.
  Inst. Ser. C: Math. Phys. Sci.}, pages 131--140. Kluwer, 1989.

\bibitem{bf91}
G.~Bj\"{o}rck and R.~Fr\"{o}berg.
\newblock A faster way to count the solutions of inhomogeneous systems of
  algebraic equations, with applications to cyclic n-roots.
\newblock {\em J. Symbolic Computation}, 12(3):329--336, 1991.

\bibitem{bf94}
G.~Bj\"{o}rck and R.~Fr\"{o}berg.
\newblock Methods to ``divide out'' certain solutions from systems of algebraic
  equations, applied to find all cyclic 8-roots.
\newblock In M.~Gyllenberg and L.E. Persson, editors, {\em Analysis, Algebra
  and Computers in Math. research}, volume 564 of {\em Lecture Notes in Applied
  Mathematics}, pages 57--70. Marcel Dekker, 1994.

\bibitem{bgk86}
W.~Boege, R.~Gebauer, and H.~Kredel.
\newblock Some examples for solving systems of algebraic equations by
  calculating {G}roebner bases.
\newblock {\em J. Symbolic Computation}, 2:83--98, 1986.

\bibitem{but84}
C.~Butcher.
\newblock An application of the {R}unge-{K}utta space.
\newblock {\em BIT}, 24:425--440, 1984.

\bibitem{cd88}
H~Cohn and J.~Deutsch.
\newblock An explicit modular equation in two variables for {Q}[sqrt(3)].
\newblock {\em Math. Comp.}, 50:557--568, 1988.

\bibitem{dav87}
J~Davenport.
\newblock Looking at a set of equations.
\newblock Technical report 87-06, Bath Computer Science, 1987.

\bibitem{emi94}
I.Z. Emiris.
\newblock {\em Sparse Elimination and Applications in Kinematics}.
\newblock PhD thesis, Computer Science Division, Dept. of Electrical
  Engineering and Computer Science, University of California, Berkeley, 1994.

\bibitem{ec95}
I.Z. Emiris and J.F. Canny.
\newblock Efficient incremental algorithms for the sparse resultant and the
  mixed volume.
\newblock {\em J. Symbolic Computation}, 20(2):117--149, 1995.
\newblock Software available at {\tt http://www.inria.fr/saga/emiris}.

\bibitem{ev97}
I.Z. Emiris and J.~Verschelde.
\newblock How to count efficiently all affine roots of a polynomial system.
\newblock Rapport de recherche no. 3212, INRIA, 1997.
\newblock To appear in {\em Discrete Applied Mathematics.}

\bibitem{Fischer}
G.~Fischer.
\newblock {\em Complex Analytic Geometry}, volume 538 of {\em Lecture Notes in
  Mathematics}.
\newblock Springer--Verlag, New York, 1976.

\bibitem{Fulton}
W.~Fulton.
\newblock {\em Intersection Theory}, volume (3) 2 of {\em Ergeb.\ Math.\
  Grenzgeb.}
\newblock Springer--Verlag, Berlin, 1984.

\bibitem{gt89}
A.~Galligo and C.~Traverso.
\newblock Practical determination of the dimension of an algebraic variety.
\newblock In E.~Kaltofen and S.M. Watt, editors, {\em Computers and
  Mathematics}, pages 46--52. Springer--Verlag, Berlin, New York, 1989.

\bibitem{glw98}
T.~Gao, T.Y. Li, and X.~Wang.
\newblock Finding isolated zeros of polynomial systems in ${C}^n$ with stable
  mixed volumes.
\newblock To appear in {\em J. of Symbolic Computation}.

\bibitem{Gunning}
R.C. Gunning.
\newblock {\em Lectures on Complex Analytic Varieties: The Local
  Parameterization Theorem}.
\newblock Princeton University Press, Princeton, 1970.

\bibitem{haa96}
U.~Haagerup.
\newblock Orthogonal {MASA}'s in the n x n - matrices, complex {H}adamard
  matrices and cyclic n-roots.
\newblock Talk at the conference on ``Operator Algebras and Quantum Field
  Theory'', Rome, July 1-6, 1996.

\bibitem{hs97}
B.~Huber and B.~Sturmfels.
\newblock Bernstein's theorem in affine space.
\newblock {\em Discrete Comput. Geom.}, 17(2):137--141, 1997.

\bibitem{li97}
T.Y. Li.
\newblock Numerical solution of multivariate polynomial systems by homotopy
  continuation methods.
\newblock {\em Acta Numerica}, 6:399--436, 1997.

\bibitem{lsy89}
T.Y. Li, T.~Sauer, and J.A. Yorke.
\newblock The cheater's homotopy: An efficient procedure for solving systems of
  polynomial equations.
\newblock {\em SIAM J. of Numerical Analysis}, 26:1241--1251, 1989.

\bibitem{lw96}
T.Y. Li and X.~Wang.
\newblock The {BKK} root count in {$C^n$}.
\newblock {\em Math. Comp.}, 65(216):1477--1484, 1996.

\bibitem{moe98}
H.M. M{\"{o}}ller.
\newblock Gr{\"{o}}bner bases and numerical analysis.
\newblock In B.~Buchberger and F.~Winkler, editors, {\em Gr{\"{o}}bner Bases
  and Applications}, volume 251 of {\em London Mathematical Lecture Note
  Series}, pages 159--178. Cambridge University Press, 1998.

\bibitem{mor87}
A.P. Morgan.
\newblock {\em Solving polynomial systems using continuation for engineering
  and scientific problems}.
\newblock Prentice-Hall, Englewood Cliffs, N.J., 1987.

\bibitem{MS1}
A.P Morgan and A.J. Sommese.
\newblock A homotopy for solving general polynomial systems that respects
  $m$-homogeneous structures.
\newblock {\em Appl. Math. Comput.}, 24:101--113, 1987.

\bibitem{MS2}
A.P Morgan and A.J. Sommese.
\newblock Computation of all solutions to polynomial systems using homotopy
  continuation.
  \newblock {\em Appl. Math. Comput.}, 24:115--138, 1987.
  \newblock Errata: {\em Appl. Math. Comput.} 51 (1992), p. 209.

\bibitem{MS3}
A.P. Morgan and A.J. Sommese.
\newblock Coefficient-parameter polynomial continuation.
\newblock {\em Appl. Math. Comput.}, 29:123--160, 1989.
\newblock Errata: {\em Appl. Math. Comput.} 51 (1992), p. 207.

\bibitem{mw90}
A.P. Morgan and C.W. Wampler.
\newblock Solving a planar four-bar design problem using continuation.
\newblock {\em ASME J. of Mechanical Design}, 112:544--550, 1990.

\bibitem{roj94}
J.M. Rojas.
\newblock A convex geometric approach to counting the roots of a polynomial
  system.
\newblock {\em Theoret. Comput. Sci.}, 133(1):105--140, 1994.
\newblock Extensions and corrections available at
  http://www.cityu.edu.hk/ma/staff/rojas/.

\bibitem{roj99}
J.M. Rojas.
\newblock Toric intersection theory for affine root counting.
\newblock {\em Journal of Pure and Applied Algebra}, 136(1):67--100, 1999.

\bibitem{rw96}
J.M. Rojas and X.~Wang.
\newblock Counting affine roots roots of polynomial systems via pointed
  {N}ewton polytopes.
\newblock {\em J. Complexity}, 12:116--133, 1996.

\bibitem{SW}
A.J. Sommese and C.W. Wampler.
\newblock Numerical algebraic geometry.
\newblock In J.~Renegar, M.~Shub, and S.~Smale, editors, {\em The Mathematics
  of Numerical Analysis}, volume~32 of {\em Lectures in Applied Mathematics},
  pages 749--763, 1996.
\newblock Proceedings of the AMS-SIAM Summer Seminar in Applied Mathematics,
  Park City, Utah, July 17-August 11, 1995, Park City, Utah.

\bibitem{ver97}
J.~Verschelde.
\newblock {PHC}pack: A general-purpose solver for polynomial systems by
  homotopy continuation.
\newblock To appear in {\em ACM Transactions on Mathematical Software}. Paper
  and software available at {\tt http://www.math.msu.edu/{\~{}}jan}.

\bibitem{wam96}
C.W. Wampler.
\newblock Isotropic coordinates, circularity and {B}ezout numbers: planar
  kinematics from a new perspective.
\newblock Proceedings of the 1996 ASME Design Engineering Technical Conference,
  Irvine, California August 18-22, 1996. CD-ROM edited by McCarthy, J.M.,
  American society of mechanical engineers. Also available as GM Technical
  Report, Publication R{\&}D-8188., 1996.

\end{thebibliography}
\end{document}